\documentclass[12pt]{book}
\title{
\textbf{From Philosophy to Program Size \\ 
Key Ideas and Methods\vspace{5mm} \\
Lecture Notes on \\
Algorithmic Information Theory\vspace{1cm} \\
Estonian Winter School in 
\\ Computer Science \\
March 2003}
}

\author{
Gregory Chaitin, \emph{IBM Research}
}

\date{}

\begin{document}
\maketitle
\markboth{CHAITIN, FROM PHILOSOPHY TO PROGRAM SIZE}{}

\chapter*{Preface}
    
This little book contains the course that I had the pleasure of giving
at the Eighth Estonian Winter School in Computer Science (EWSCS '03)
held at the 
beautiful
Park Hotel Palmse in Lahemaa National Park, Estonia, from March 2nd through 7th, 2003.
There I gave
four 90-minute lectures 
on algorithmic information theory (AIT), which is 
the theory of program-size complexity.
Each of these lectures is one chapter of this book.
     
In these lectures I discuss philosophical applications of AIT,
not practical applications.
Indeed, I believe AIT has no practical applications.
    
The most interesting thing about AIT is that you can almost never determine the
complexity of anything. This makes the theory useless for practical
applications, but fascinating from a philosophical point of view, because it
shows that there are limits to knowledge.

Most work on computational complexity is concerned with time. However this
course will try to show that program-size complexity, which measures
algorithmic information, is of much greater philosophical significance. 
I'll discuss how one can use this complexity measure to study what
can and cannot be achieved by formal axiomatic mathematical theories. 

In
particular, I'll show \textbf{(a)} that there are natural information-theoretic
constraints on formal axiomatic theories, and that program-size complexity
provides an alternative path to incompleteness from the one originally used by
Kurt G\"odel. Furthermore, I'll show \textbf{(b)} that in pure mathematics there are
mathematical facts that are true for no reason, that are true by accident.
These have to do with determining the successive binary digits of the precise
numerical value of the halting probability $\Omega$ for a ``self-delimiting'' universal
Turing machine. 

I believe that these meta-theorems \textbf{(a,b)} showing \textbf{(a)} that the
complexity of axiomatic theories can be characterized information-theoretically
and \textbf{(b)} that God plays dice in pure mathematics, both strongly suggest a
quasi-empirical view of mathematics. I.e., math is different from physics,
but perhaps not as different as people usually think. 

I'll also discuss the
convergence of theoretical computer science with theoretical physics, Leibniz's
ideas on complexity, Stephen Wolfram's book \emph{A New Kind of Science}, and
how to attempt to use information theory to define what a living being is.

In this book I've tried to preserve the informal style of presentation in my lectures
that stressed the key ideas and methods and avoided getting bogged down in
technical details. There are no proofs here, but there are plenty of proof sketches.
I hope that you enjoy reading this book just as much as I enjoyed presenting 
this material at EWSCS '03!
\vspace{5mm}
\\
\emph{---Gregory Chaitin}
\vspace{2mm}
\\
\texttt{http://www.cs.auckland.ac.nz/CDMTCS/chaitin}

\tableofcontents

\chapter{Day I---Philosophical Necessity of AIT}
\markboth{CHAITIN, FROM PHILOSOPHY TO PROGRAM SIZE}{DAY I---PHILOSOPHY}

\section{Is $\mathbf{P} \neq \mathbf{NP}$ a New Axiom?}
    
Can we add $\mathbf{P} \neq \mathbf{NP}$ as a new axiom?!
    
This is a good example of the situation discussed in 
G\"odel, \emph{What is Cantor's Continuum Problem?}, 1947,
where he argues that \textbf{maybe math is a little like physics} and
that new axioms that are not self-evident might be justified 
because of their \textbf{usefulness}.
   
If so, there is ample pragmatic justification for adding
$\mathbf{P} \neq \mathbf{NP}$
as a new axiom.
   
(In his 1947 article G\"odel was concerned with set theory, 
not computer science.)
    
Let's see!!!

\section{Day I Summary}
    
\begin{itemize}
\item
Our goal: To find the
\emph{Limits of Mathematical Reasoning}
\item
Tool we'll use:
\emph{Program-Size \textbf{Complexity}}
\item
It's also called:
\emph{Algorithmic \textbf{Information}}
\item
Idea starts with \textbf{Leibniz!} 1686! 
\\   
\emph{Discourse on Metaphysics} (in French!)
\item
AIT is needed to understand what is ``law'' and 
\\
what is ``understanding.''
\end{itemize}

\section{Summary of Leibniz, 1686}
    
(Brought to my attention by Hermann 
Weyl, 1932.\footnote
{Hermann Weyl, \emph{The Open World}, Yale University Press, 1932. Reprinted by Ox Bow Press, 1989.
See also Weyl's \emph{Philosophy of Mathematics and Natural Science}, Princeton University Press, 
1949.})
    
What is a law of nature?
    
According to Leibniz,
a theory must be \textbf{simpler} 
than the data it explains!
    
Because \textbf{if a physical law can be as
complicated as the experimental data
that it explains, then
there is always a law, and the
notion of ``law'' becomes \emph{meaningless}!}
    
Understanding is compression!
A theory as \textbf{complicated} as the data it explains is \textbf{NO} theory!
    
All of this is stated very clearly (in French) in 1686
by the mathematical genius and philosopher Leibniz!
During his lifetime he only transmitted summaries
of these ideas to friends and colleagues in letters!
The complete text of his \emph{Discourse on Metaphysics}
was found among his voluminous personal papers after his death.\footnote
{For the original French text, see Leibniz, \emph{Discours de m\'etaphysique}, Gallimard, 1995. 
There is an English translation in G. W. Leibniz, \emph{Philosophical Essays}, 
edited and translated by Roger Ariew and Daniel Garber, Hackett,
1989.}

\section{From Leibniz to AIT}
    
AIT goes beyond Leibniz by positing 
that \textbf{a theory is a computer program},
and that \textbf{the size in bits of this computer
program is the complexity of the theory}.
AIT makes Leibniz more precise!
AIT emphasizes the common features in these four diagrams:
\begin{center}
scientific theory $\longrightarrow$ \textbf{Calculations} $\longrightarrow$ experimental data
\vspace{1mm}
\\
program $\longrightarrow$ \textbf{Computer} $\longrightarrow$ output
\vspace{1mm}
\\
axioms $\longrightarrow$ \textbf{Deduction} $\longrightarrow$ theorems
\vspace{1mm}
\\
Ideas $\longrightarrow$ \textbf{Mind of God} $\longrightarrow$ The World
\end{center}
    
\textbf{Following Leibniz, in each of these diagrams
the input on the left 
must be much smaller 
than the output on the right.}
    
Leibniz's key insight is not 
that this is ``the best of all possible worlds''.
This was anti-Leibniz propaganda by Voltaire, who ridiculed Leibniz and
did not understand how subtle and profound Leibniz was.
(According to Borges, the word ``optimism'' was invented by Voltaire to mock Leibniz!)
     
Leibniz's key insight 
is that God has used few ideas to create all
the diversity, richness and apparent complexity of the natural world.\footnote
{Here is an argument in favor of diversity from Leibniz via Borges.
Consider two libraries with exactly the same number of books. Recall that Borges was
a librarian!  One of the libraries only has copies of Virgil's \emph{Aeneid} (presumably
the perfect book), while the other library has \textbf{one} copy of the \emph{Aeneid} and
many other different (and presumably inferior!)\@ books. Nevertheless, the second 
library is clearly a more interesting place to visit!}
Leibniz is actually affirming his belief that the universe is \textbf{rationally comprehensible}.
(This belief in a rational universal goes back at least to the ancient Greeks,
particularly Pythagoras and Plato, but Leibniz's formulation is much sharper
and profound because he analyzes in mathematical terms exactly what this belief \textbf{means}.)
In modern language, Leibniz was stating \textbf{his belief in the possibility of science}.
     
Pythagoras and Plato believed that the universe can be comprehended using mathematics.
Leibniz went beyond them by clarifying mathematically what exactly does it mean
to assert that the universe can be comprehended using mathematics.
    
AIT continues this train of thought and goes beyond Leibniz by positing that explanations
are computer programs and also by \textbf{defining} 
precisely what \textbf{complexity} is and what exactly
does it mean to satisfy Leibniz's requirement that 
\textbf{an explanation has to be simpler than
the phenomena that it explains}. 

\section{Setup for Discussing the Limits of Mathematical Reasoning}
    
Now let's \textbf{apply} these ideas!
Let's see what they have to say about the limits of mathematical reasoning.
Let's see what mathematical theories can accomplish.
    
The setup is as follows.
    
The (static, formal) mathematical theories that we consider
are thought of, somewhat abstractly, as a computer 
program for running through all possible proofs and generating
all the theorems in the theory, which are precisely all the
logical consequences of the axioms in the theory.
\begin{center}
Theory = Program  
$\longrightarrow$ 
\textbf{Computer}
$\longrightarrow$ 
Theorems
\end{center}
    
[This assumes that the theory uses a formal language, symbolic logic,
requires complete proofs, and provides a proof-checking algorithm
that always enables you to decide mechanically whether or not a proof
is correct, whether or not it followed all the rules in deriving 
a logical consequence of the axioms.]
    
And we shall concentrate our attention on the \textbf{size in bits} of this computer
program for generating all the theorems.  That's our measure of the \textbf{complexity}
of this theory, that's our measure of its \textbf{information content}.
    
So when we speak of an $N$-bit theory, we mean one with an $N$-bit program
for generating all the theorems.  We don't care that this process is very
slow and never terminates. AIT is a theoretical theory, not a practical theory!
    
Okay, that's half the setup!  These are the theories we'll consider.
\vspace{1cm}
    
Here's the other half of the setup.  
Here are the theorems that we want these theories 
to be able to establish.
We want to use these theories to prove that 
particular computer programs are \textbf{elegant}.
\begin{center}
``Elegant'' Programs:
\\
\emph{No smaller program gives exactly the same output.}
\end{center}

What is an elegant program?  It's one with the property
that no smaller program written in the same language
produces exactly the same output.\footnote
{Throughout this discussion we assume a \textbf{fixed choice}
of programming language.}
    
There are lots of elegant programs!
    
For any particular computational task, for any particular
output that we desire to achieve, there has to be at 
least \textbf{one} elegant program, and there may even be several.  

But what if we want to be able to \textbf{prove} that a particular
program is elegant?

\section{You need an $N$-bit theory 
to prove that an $N$-bit program 
is elegant!}
    
[For the proof, see Section 1.8; for a corollary, see Section 1.7.]
    
These are
\textbf{\emph{Irreducible} Mathematical Truths}!!!
    
No rational justification is possible!
    
Such mathematical facts can have no 
rational explanation, because rational 
explanation consists of reducing something
to simpler (maybe even self-evident)
principles.  A theory for something
derives it from \textbf{simpler} hypotheses.
If this kind of reduction is impossible,
then we are confronted with \textbf{irrational}
mathematical facts, mathematical facts
that cannot be encompassed by \textbf{any} static theory!

Therefore (corollary)\ldots

\section{Complexity of the Universe of Mathematical Ideas}
    
The world of mathematical ideas
has \textbf{INFINITE} complexity! 
    
Why?
    
Well, no $N$-bit theory for any \textbf{finite} $N$ can enable
you to prove all true assertions of the form
\begin{center}
``This particular program is elegant.''
\end{center}
    
(What about the \textbf{physical} world?
Does it have finite or infinite complexity?
We'll look at that later.
See Sections 1.10 through 1.12.)

\section{Why do you need an $N$-bit theory to prove 
that an $N$-bit program is elegant?}
    
Here is the proof!
We'll assume the opposite and derive a contradiction.
Consider this computer program:
\newpage
\begin{center}     
$N$-bit theory + fixed-size routine
\\   
$\downarrow$ 
\\  
\textbf{Computer} 
\\   
$\downarrow$ 
\\  
\emph{the output of}
\\
the first provably elegant program
whose size is greater than 
\\
(complexity of theory + 
size of fixed-size routine)
\end{center}
     
Let $c$ be the size in bits of the fixed-size routine
that does all this: It is given the $N$-bit theory as data and 
then it runs through all possible proofs in
the theory until it finds a provably elegant program that
is sufficiently large ($> N + c$ bits), then it runs that program and returns
that large elegant program's output as its own output.
   
So the $N + c$ bit program displayed above produces 
the same output as a provably elegant program that is larger
than $N + c$ bits.  But that is impossible! 
    
So our precise result is that \textbf{an $N$-bit theory
cannot enable us to prove that any program
that is larger than $N + c$ bits is elegant}.
Here $c$ is the size in bits of the fixed-size
routine that when added to the $N$-bit formal theory
as above yields the paradox that proves our theorem.
    
\textbf{Note:} We are making the tacit assumption that
if our theory proves that a program is elegant, then that
program is actually elegant.  I.e., we assume that only
true theorems are provable in our formal theory.  If this is
\textbf{NOT} the case, then the theory is of absolutely no interest.

\section{Is Mathematics Quasi-Empirical?}
    
Now let's go back to the question of adding 
$\mathbf{P} \neq \mathbf{NP}$
as a new axiom,
and to G\"odel's thoughts that maybe
physics and math are not so different,
which, following Lakatos and Tymoczko (1998),
is now referred to as a quasi-empirical
view of mathematics.
    
G\"odel, who had the conventional Platonist view of math,
was only forced into backing \textbf{new
math axioms that are only justified pragmatically}
just as in physics because of his famous 1931 incompleteness
theorem.  And I believe that the ideas that
I've just presented applying program-size
complexity to incompleteness, in particular
my result that it takes an $N$-bit theory
to prove that an $N$-bit program is elegant,
and the results on $\Omega$ that we'll see Day II,
provide even more support for
\textbf{G\"odel's heretical views on new axioms}.
     
The way AIT measures the complexity (information
content) of mathematical theories makes G\"odel
incompleteness seem much more natural, much more
pervasive, much more inevitable, and much more
\textbf{dangerous}. Adding new axioms---adding
more mathematical information---seems to be
the only way out, the only way to go forward!  
    
We've discussed adding new axioms in math just
as in physics, pragmatically. 
A related question is 
\begin{center}
``Is \textbf{experimental mathematics} okay?''
\end{center}
Even when there are \textbf{NO} proofs?
For an \textbf{extreme example} of this, see
Wolfram, \emph{A New Kind of 
Science}, 2002,
who provides a tremendous amount of 
computational evidence, but almost no proofs,
in support of his theory.
See also the journal \emph{Experimental Mathematics}.
   
Obviously, AIT makes me sympathetic to experimental mathematics,
even though I don't do experimental math myself.  Experimental math is fueled
by the power of the computer, not by G\"odel's nor by my 
meta-theorems, but it fits nicely into a quasi-empirical view
of mathematics. Practical necessity as well as philosophical
analysis are both simultaneously pushing in the direction
of experimental math!

\section{Complexity of the Physical World}
    
Now let's turn to the \textbf{physical} universe!
Does it have finite or infinite complexity?
    
The conventional view on this held by high-energy physicists
is that there is a TOE, a theory of everything, a finite set
of laws of nature that we may someday know, which has
only finite complexity. 
    
So that part is optimistic!
    
But unfortunately in quantum mechanics there's
randomness, God plays dice, and to know the results of all 
God's coin tosses, infinitely many coin tosses,
necessitates a theory of infinite complexity, which 
simply records the result of each toss!
    
So the conventional view currently held by physicists
is that because of randomness in quantum mechanics the world has infinite complexity.
    
Could the conventional view be wrong?
Might it nevertheless be the case that the universe only has finite complexity?
Some extremely interesting thoughts on this can
be found in
\begin{center}
\textbf{Wolfram, \emph{A New Kind of Science}, 2002.}
\end{center}
    
According to Wolfram, 
there is \textbf{no real randomness}.
There is only \textbf{pseudo-randomness},
like the randomness produced by random-number
generators, which are actually deterministic
sequences of numbers governed by mathematical
laws, since computers use algorithms to generate
pseudo-random numbers, they don't use quantum mechanics!
    
According to Wolfram,
our universe is actually a \textbf{deterministic universe}
that's governed by deterministic physical laws!
    
So the physical universe, the world, has \textbf{finite complexity}.
According to Wolfram,
\textbf{everything happens for a reason},
just as Leibniz thought! 
    
These two supremely intelligent men are rationalists.  
They want to understand \textbf{everything}!
They don't believe in ultimate mysteries!
They don't think anything is incomprehensible!
They believe in the power of the human mind to
comprehend everything!
    
\textbf{On the other hand, we have seen that 
because it has infinite complexity,
the universe of mathematical ideas
CANNOT be comprehended in its entirety.}

\section{Summary of Wolfram, 2002}

In summary, according to
Wolfram, 2002
the world is like 
$\pi$ = 3.1415926\ldots
    
It \textbf{looks complicated},\footnote 
{Especially if you're given successive digits 
of $\pi$ that come from far inside the decimal
expansion without being told that they're digits of $\pi$---they look random.}
but it is actually 
very \textbf{simple}!
    
According to Wolfram
all the randomness we see
in the physical world is
actually \textbf{pseudo}-randomness.
He believes that the physical world
is actually deterministic,
we just don't know the law.
    
He sees this as a \textbf{philosophical possibility}.
Whether \textbf{our} physical universe is \textbf{actually}
that way or not is another matter,
to be decided by scientists, not philosophers!

\section{Digital Philosophy}
    
Wolfram's book as well as my 
own work on AIT are both examples of what
Edward Fredkin refers to as
\textbf{digital philosophy}, a viewpoint that Fredkin also helped to pioneer.\footnote
{For more on Wolfram, Fredkin, Lloyd, Toffoli, Landauer, Zuse\ldots\ see
O. Postel-Vinay, ``L'Univers est-il un calculateur?''\@ [Is the universe a calculator?], 
\emph{La Recherche}, no.\ 360, January 2003, pp.\ 33--44.} 
    
In a nutshell, \emph{digital philosophy} posits that
the world is a giant computer,
a giant digital information processor,
and that, fundamentally, 
\textbf{everything is discrete 0/1 bits}!
     
This algorithmic view of everything\footnote 
{Algorithms played a decisive role in Sumerian mathematics
more than a \textbf{millennium before} Pythagoras, a tremendously long intellectual trajectory!
The Sumerians used base
60 numerals, and divided the circle into 360 degrees.\emph{---Fran\c{c}oise 
Chaitin-Chatelin, private communication.}}
works much better if there are
actually no real numbers, no continuous quantities, and the physical
universe is really, at some bottom level, discrete.
    
Wolfram's work, AIT, and Fredkin's digital philosophy are all examples of
the \textbf{convergence of mathematics, 
theoretical physics, and theoretical computer science}!
This is an accelerating trend, of which the field of quantum computing
is also an example.
    
Of course, traditional mathematical physics is based 
on continuous math,
on ordinary and partial differential equations,
and does not fit in too well with a digital philosophy.
Maybe digital philosophy is a terrible mistake.
Maybe we are taking the digital computer much too seriously!
Maybe we shouldn't make it the basis of a new philosophy, of a new world view, of
a new syst\`eme du monde?
    
We will see!

\section{Two Interesting Books}
    
I should mention that besides my own works,
\textbf{the} book on the quasi-empirical view of math is
\begin{center}
\textbf{
Tymoczko,
\\
\emph{New Directions in the 
Philosophy of Mathematics,}
\\
Princeton University Press,
1998.
}
\end{center}
     
This is the second expanded edition of a valuable collection
of essays by philosophers, mathematicians, and computer
scientists (including two of my own essays) that its editor
Tymoczko unfortunately did not live to see in print.
    
Highly recommended!
    
There is also a forthcoming book on experimental math,
\begin{center}
\textbf{
Bailey (and Borwein?), 
\\
\emph{The Experimental Mathematician}, 
\\A. K. Peters, 2003?,
}
\end{center}
that should be extremely interesting.

\section{Decomposing the World}
    
Finally, here is a new and completely different philosophical application of AIT.
\vspace{5mm}
    
What is a living being?
    
How can we partition the world into parts?
Can we do this
in spite of Parmenides and mystics who insist
that the world must be apprehended as a whole
and is an organic unity, a single substance,
and \textbf{cannot} be separated into independent parts?
\vspace{5mm}
    
I think the key is
\textbf{algorithmic independence}.
     
$X$ and $Y$ are said to be
algorithmically independent if 
the program-size complexity of $X$ \textbf{and} $Y$
$\approx$ the sum of the individual program-size complexities
of $X$ and $Y$.
    
I.e., if the number of bits in the simplest theory that explains both simultaneously
is approximately equal to the sum of the number of bits in the simplest theories that explain
each of them separately.
\vspace{5mm}
     
Independent parts of the world, of which living beings are the most interesting example,
have the property that their program-size complexity 
\textbf{decomposes additively}.
    
Conversely, the parts of a living being are \textbf{certainly not independent} and
have high \textbf{mutual information}.  [Mutual information will be discussed
in Sections 3.6 and 3.7, Day III, and Section 4.3, Day IV.]
\vspace{5mm}
    
\emph{This needs much more work!}

\chapter{Day II---Main Application of AIT: Incompleteness}
\markboth{CHAITIN, FROM PHILOSOPHY TO PROGRAM SIZE}{DAY II---INCOMPLETENESS}

\section{Day II Summary}
    
Limits of Formal 
Mathematical Reasoning:
\begin{itemize}
\item[\emph{I.}]
You can't prove a number is 
uninteresting/random.
\item[\emph{II.}]
You can't determine bits of $\Omega$
(accidental mathematical facts).
\end{itemize}
Both \emph{I} and \emph{II} are 
\textbf{\emph{irreducible}} mathematical truths!

\section{Hilbert-Style Formal Theories}

\begin{itemize}
\item
\textbf{Meta-mathematics:}
Use mathematics to discover the limits of
mathematics.
\textbf{Cannot} be applied to intuitive, informal mathematics.

\item
What's a \textbf{formal theory}?
\begin{itemize}
\item
axioms
\item
rules of inference
\item
symbolic logic
\item
formal grammar
\end{itemize}
There is a proof-checking algorithm!
    
\item
There is algorithm for generating all theorems 
in size order of proofs! So, following Emil Post, 1944,
formal theory = r.e.\ set of propositions:
\begin{center}
theory $\longrightarrow$ \textbf{Computer} $\longrightarrow$ theorems
\end{center}
   
\item
This is a \textbf{\emph{static}} view of mathematics!
It is assumed that the axioms and rules of inference don't change.
Mathematical methods \textbf{can be used} to
discuss what such a formal system can achieve.
\end{itemize}

\section{Summary of G\"odel, 1931}
    
\begin{itemize}
\item
``This statement is false!''
    
\item
True iff false! Paradox!

\item
``This statement is unprovable${}_{FT}$!''

\item
True iff unprovable${}_{FT}$!
Therefore either a false statement is provable${}_{FT}$, and \emph{FT}
is useless, or a true statement is unprovable${}_{FT}$, and \emph{FT} is incomplete!

\item
\emph{FT} = 
\emph{F}ormal 
\emph{T}heory. 
\end{itemize}

\section{My Approach to Incompleteness}
    
\begin{itemize}
\item
Let's assume we can
divide all the positive integers into the interesting and the uninteresting ones,
and that infinitely many of them are uninteresting.

\item
Consider
the first uninteresting positive integer!
   
\item
It's
\emph{ipso facto} interesting!
Paradox!
    
\item
Consider
the first provably${}_{FT}$ uninteresting positive integer!
    
\item
It cannot exist because it would be an extremely interesting number!
Therefore \emph{FT} is incomplete!
\emph{FT} can \textbf{never} prove that a number is uninteresting!
[Actually, it will turn out, \textbf{almost} never.]
    
\item{}
How do we actually define an uninteresting number?
\emph{Uninteresting} 
means algorithmically irreducible,
incompressible, no theory for it smaller than it is.\footnote
{I.e., using notation that we haven't introduced yet,
the positive integer $N$ is uninteresting iff $H(N) \equiv H_U(N) \geq |N| \equiv$
size in bits of $N$.
For the definition of $H_U$, see Section 2.13.}
     
\item
\textbf{Final Conclusion:}
Infinitely many integers are uninteresting.
But using a \textbf{fixed} \emph{FT} you can almost never prove it!
Only in finitely many cases!\footnote{
It turns out that if you work through all the details,
you can't prove that a positive integer is uninteresting 
if its size in bits is \textbf{greater} than (the size
of the program that generates all the theorems in
the formal theory \emph{FT}) plus a constant $c$.}
\end{itemize}

\section{Borel's Unreal Number $B$, 1927}

Let's start on the path to
the halting probability $\Omega$, which is a
real number.

\begin{itemize}
\item
Start with
Borel's \textbf{very unreal} real number $B$.
    
\item
The original source is \'Emile Borel, 1927.\footnote 
{In English in Paolo Mancosu, \emph{From Brouwer to Hilbert}, Oxford University Press, 1998.}
\\
$B$ was brought to my attention by Vladimir Tasi\'c, 2001.\footnote
{Vladimir Tasi\'c, \emph{Mathematics and the Roots of Postmodern Thought}, 
Oxford University Press, 2001.}
     
\item
$B = 0. d_1 d_2 d_3 \ldots$
     
\item
The $N$th digit $d_N$ of $B$ answers the $N$th question in French!
     
\item
$d_N =$
\begin{itemize}
\item[$1 \longrightarrow$]
answer \emph{\textbf{Yes}!}
\item[$2 \longrightarrow$]
answer \emph{\textbf{No}!}
\item[$0 \longrightarrow$]
\emph{\textbf{NO ANSWER}}
\\
($N$th text in French 
is not a valid yes/no question, 
\\
or cannot be answered.)
\end{itemize}
    
\item
Most digits of Borel's number $B$ are 0's.
$B$ is extremely redundant!
\end{itemize}

\section{More-Real Turing Number $T$}

Here is the next step on our path to $\Omega$:

\begin{itemize}
\item
\textbf{Real Number $T$ Solving Turing's Halting Problem}
    
\item
$T = 0. b_1 b_2 b_3 \ldots$
     
\item
The $N$th bit $b_N$  of $T$
answers the $N$th case of the halting problem.
\\   
$b_N$ tells us whether the $N$th program $p_N$ 
ever halts.
     
\item
$b_N = $  
\begin{itemize}
\item[$0 \longrightarrow$]
$p_N$ doesn't halt,
\item[$1 \longrightarrow$]
$p_N$ halts.
\end{itemize}
\end{itemize}

\section{Some Interesting Cases of the Halting Problem}
    
\textbf{These are bits of $T$:}
   
\begin{itemize}
\item
\emph{Fermat's Last Theorem} (Andrew Wiles)
\\   
Does 
\[
x^N + y^N = z^N
\]
have a solution 
in positive integers with $N \ge 3$?
\item
\emph{Riemann Hypothesis}
\\   
About the location of the complex zeroes of the zeta function
\[   
\zeta(s)\equiv \sum_n \frac{1}{n^s} 
= \prod_p \frac{1}{1-\frac{1}{p^s}}
\]   
(Here $n$ ranges over positive integers and $p$ ranges over the primes.)
\\  
Tells us a lot about the distribution of prime numbers.
\item
\emph{Goldbach Conjecture}
\\   
Is every even number the sum of two primes?
\end{itemize}
\textbf{In each case there is a program that systematically searches
for a counter-example and halts iff it finds it.}

\section{Crucial Observation}
    
Suppose we are given $N$ programs
and want to know 
which ones halt and
which ones don't.
     
\textbf{$N$ cases of the halting problem is
only $\log_2 N$ bits of information,
not $N$ bits!}
    
They are \textbf{never} independent mathematical facts!
    
Why not?
    
Because
we could answer all $N$ cases of the halting problem
if we knew \textbf{exactly how many} of the $N$ programs halt.
Just run them all in parallel until exactly the right number halt.
The others will \textbf{never} halt.
    
Now we use this observation to compress all the redundancy
out of the Turing number $T$ and get an algorithmically
irreducible number $\Omega$.

\section{The Halting Probability $\Omega$}

Finally $\Omega$ = Omega!
    
\textbf{Halting Probability:}
    
\[
\Omega \equiv \sum_{\mbox{$p$ halts}} 2^{-|p|}
\]
    
$|p|$ is the size in bits of the program $p$.
    
I.e., each $k$-bit program that halts
when run 
on our standard universal Turing machine $U$
contributes $1/2^k$ to $\Omega$.
    
[We need to make $U$ \textbf{self-delimiting} (Section 2.12)
to ensure that $\Omega \le 1$. Otherwise the sum
for $\Omega$ diverges to $\infty$.
By using self-delimiting 
programs, we've constructed \textbf{one} number, $\Omega$,
that extends the trick of Section 2.8's \emph{crucial observation} 
so that it works for
an \textbf{infinite} number of 
computer programs, all of them, in fact.]
   
Now there is absolutely \textbf{NO} redundancy! 
   
The first $N$ bits of $\Omega$
answers the halting problem for all
programs up to $N$ bits in size!
(Can you see why? \emph{Hint:} $\Omega$ can be computed in
the limit from below.\footnote
{But \textbf{very, very} slowly, and you can never be sure how close you are.})

\section{Why is $\Omega$ Interesting?}

The base-two bits of $\Omega$ are 
\textbf{irreducible mathematical facts!}
They can't be derived from anything simpler!
    
The bits of $\Omega$ are algorithmically
irreducible, 
algorithmically independent,
and
algorithmically random!
    
$\Omega$ is a real number with maximal information
content.  
Each bit of $\Omega$ is a complete surprise!
It's not at all a \textbf{tame} real number like $\pi$ = 3.1415926\ldots\
which only has a finite amount of algorithmic information.
    
$\Omega$ is a \textbf{dangerous, scary} real number!
Not good to meet in a dark alley!
$\Omega$ is maximally unknowable!
Maximum entropy, maximum disorder!
\begin{center}
``\textbf{$\Omega$ is a nightmare for the rational mind!}'', \emph{Karl Svozil}.
\end{center}
        
This makes $\Omega$ sound bad, very bad!
    
On the other hand, $\Omega$ is distilled, crystalized mathematical
information.  
If $T$ is coal, then $\Omega$ is a diamond!

In fact, initial segments of $\Omega$ are ideal \textbf{new mathematical axioms}.
Knowing a large initial segment of $\Omega$ would settle all halting problems
for programs up to that size, 
which would, to use Charles Bennett's terminology,
settle all \textbf{finitely refutable} math conjectures up to that complexity.

$\Omega$ is concentrated essence of \textbf{mathematical creativity}
and \textbf{mathematical inspiration}!
One could measure the progress of mathematics by how many bits of $\Omega$
we can currently determine!\footnote{Somewhere Leibniz proposes measuring the
intellectual progress of mankind via a function $\Phi(t)$  
with the property that all
interesting theorems with proofs of size $\leq \Phi(t)$ are known at time $t$. 
Yet another possible  
measure of human progress is a function $\Lambda(t)$ such that all halting problems for
programs with $\leq \Lambda(t)$ bits have been settled by time $t$.
Yet perhaps these measures of intellectual progress are all beside the point,
the point being the \textbf{emergence of new concepts}?

The progress measured by $\Phi$ is, in principle, merely hard work, and could be
achieved mechanically, by employing vast amounts of computation, assuming that we
are exploring a fixed, static formal theory. But I think that $\Lambda$---and counting provable
bits in $\Omega$---measures the emergence of new concepts \textbf{indirectly} via
their \textbf{effects}, for surely new concepts would be needed to advance in these areas.}
(Of course this has nothing to do with \textbf{moral} or \textbf{scientific} progress.)
    
So $\Omega$ can also be regarded as a good friend, instead of an enemy!

\section{In What Sense is $\Omega$ Random?}

\textbf{The bits of $\Omega$
are mathematical facts
that are true
for \emph{No Reason},
they're true
by \emph{Accident}!}
    
Here mathematical truth 
isn't \textbf{Black} or \textbf{White}.
It's \textbf{Grey}!
    
The best way to think about $\Omega$
is that each bit of $\Omega$ is 0 or 1
with probability 1/2!
    
Here God plays dice with 
mathematical truth!
    
Knowing all the \textbf{even} bits of $\Omega$
wouldn't help us to get any of the \textbf{odd} bits of $\Omega$!
Knowing the first \textbf{million} bits of $\Omega$
wouldn't help us to get the \textbf{million and first} bit of $\Omega$!
The bits of $\Omega$ are
just like independent tosses of a fair coin!
    
This is the case even though $\Omega$ is
a \textbf{specific} real number!
And even though each bit 
of $\Omega$
is fully determined mathematically.
(Even though we can't compute it 
nor
prove what its value is!)

\section{The Self-Delimiting Universal Computer $U$}

Actually, the value of $\Omega$ depends
on the choice of universal computer $U$.
There are many possible choices!
    
First of all, $U$ must read one bit
of its program $p$ at a time and \textbf{decide by itself}
when to stop reading $p$ before encountering
a blank at the end of $p$.
In other words, each $p$ for $U$ must be
\textbf{self-delimiting}.
    
Also, $U$ must be \textbf{universal}.
This means that for any special-purpose
self-delimiting computer $C$ there is
a prefix $\pi_C$ such that
concatenating it in front of a program for $C$
gives you a program for $U$ 
that computes the same output:
\[  
U(\pi_C \, p) = C(p).
\]   
This prefix depends only on the choice of $C$
and not on $p$.
In other words, $U$ can \textbf{simulate} each $C$.

\section{$U$ is Optimal and Defines $\Omega_U$}

Then we define the \textbf{complexity}
$H$
with respect to $C$ and $U$
as follows:
\[   
H_C(x) \equiv
\min_{C(p)=x} |p|,
\]
\[
H_U(x) \equiv
\min_{U(p)=x} |p|.
\]
$H(x)$ is the size in bits of the smallest
program 
for computing $x$ on each machine.

Then we have
\[  
H_U(x) \leq
H_C(x) + |\pi_C|.
\]
In other words, programs for $U$ are \textbf{not too large}.
For $U$ to simulate $C$ adds only a \textbf{fixed number of 
bits} to each program for $C$.

Then
$\Omega_U$ 
is defined as follows:
\[   
\Omega_U \equiv
\sum_{\mbox{$U(p)$ halts}} 2^{-|p|}.
\]
Any such universal $U$ will do.

\section{My $U$ is Programmed in LISP}

The particular $U$ that I picked to define $\Omega$
uses LISP as follows:
\begin{center}  
$U(p)$ runs $U$ on the binary program $p=\pi\beta$.
\end{center}  
Here $p = \pi\beta$ is a bit string consisting of 
a high-level algorithm
$\pi$ followed by data $\beta$.
The self-delimiting prefix $\pi$ is a LISP expression.
The data $\beta$ is raw binary data 01101\ldots\ 
The value of the LISP prefix $\pi$ is the output of $U(p)$.
    
I invented a \textbf{special version} of LISP for 
writing the
prefix $\pi$.

\section{Programming $U$ in LISP}

The LISP expression
$\pi$ is converted to binary, 
eight bits per character,
and 
concatenated with the raw binary data $\beta$
to produce the binary program for $U$, $p = \pi\beta$.
     
Access to the data $\beta$ is strictly controlled.
$U$ reads one bit at a time of $p$
and 
\textbf{CANNOT} run off the end of $p = \pi\beta$.  
The binary data $\beta$ must be
self-delimiting, 
just like the LISP prefix $\pi$, 
so that
$U$ knows when to stop reading it.
    
I.e., the alphabet for $p$ is binary, \textbf{not trinary}!
There is no blank at end of $p$!
That would be a wasted character, one that 
isn't being used nearly enough!
     
There are more details about the LISP implementation in
the last lecture (Day IV).

\section{$\Omega$ and Hilbert's 10th Problem}

We end today's lecture with a very important application:
\textbf{Hilbert's 10th problem}.
    
Does the diophantine equation 
$D(\mathbf{x}) = 0$ 
have an unsigned integer solution?
\begin{itemize}
\item
$\exists D$ such that
$D(k,\mathbf{x}) = 0$ 
has a solution iff $p_k$ halts!
\\{}
[Matijasevi\v{c}]
\item
$\exists D$ such that
$D(k,\mathbf{x}) = 0$ 
has infinitely many solutions 
\\
iff the $k$th bit of $\Omega$ is a 1!
\\{}
[Chaitin]
\item
$\exists D$ such that
$D(k,\mathbf{x}) = 0$ 
has an even/odd number of solutions 
\\
iff the $k$th bit of $\Omega$ is 0/1!
\\{}
[Ord and Kieu]
\end{itemize}
For more details, see Ord and Kieu, 
\emph{On the existence of a new family of Diophantine equations for $\Omega$},
\texttt{http://arxiv.org/abs/math.NT/0301274}.
    
\textbf{So we get \emph{randomness in arithmetic!}
$\Omega$'s algorithmic randomness also \emph{infects} 
elementary number theory!
The \emph{disease} is spreading!}

\chapter{Day III---Technical Survey of AIT: Definitions \& Theorems}
\markboth{CHAITIN, FROM PHILOSOPHY TO PROGRAM SIZE}{DAY III---DEFINITIONS \& THEOREMS}

$\Omega$ is the jewel of AIT.  But it isn't a diamond solitaire.
It's in a beautiful setting, which we'll outline today.
We'll review the key definitions, theorems, 
methods and ideas of AIT, but there will be \textbf{NO} proofs!
Only proof sketches.
For the proofs, see my book
\emph{Exploring Randomness}.

\section{Definition of a Self-Delimiting Computer $C$}

AIT starts with a
\vspace{5mm}
\\ 
\textbf{Self-delimiting computer $C$}
\begin{itemize}
\item
$C(p) \longrightarrow$ output
\item
The program $p$ is a bit string.
\item
The output can be e.g.\ a LISP S-expression \\ or a set of S-expressions.
\end{itemize}
\textbf{Abstract version of self-delimiting feature:}
\\
No extension of a valid program is a valid program.
If $C(p)$ is defined then $C$ 
is not defined on any extension of $p$.
I.e., the domain of the function $C$ is a so-called \textbf{prefix-free set}.
That's a set of words in which no word is a prefix of another.

\section{What do we do with $C$? Define $U$!}
    
First define the \textbf{program-size complexity}
$H_C(x)$ to be
the size in bits of the smallest program for $C$ to compute $x$:
\[ 
H_C(x) \equiv \min_{C(p)=x} |p|.
\]  

Next we construct a \textbf{universal computer} $U$, for example, as follows:
\[  
U(0^{\mbox{(G\"odel number for $C$)}}1p) = C(p).
\]  
In other words, $U$ can simulate $C$ if to each program for $C$ we add 
a prefix consisting of a long run of 0's followed by a 1.  
This prefix is self-delimiting, and the number of 0's in it 
indicates \textbf{which $C$ is to be simulated}.
   
Then we have
\[ 
H_U \leq H_C + \mbox{constant} 
\]  
where
\[
\mbox{constant} = 1 + \mbox{(G\"odel number for $C$)}.
\]
I.e., $U$'s programs are, to within an additive constant, minimal in size.

We take \textbf{this minimality property} as our definition for a universal computer $U$.
    
Now \textbf{somehow pick a particular natural $U$} to use as the basis for AIT.

In the lecture for Day IV, I will indicate how to \textbf{use LISP} to implement the
particular $U$ that I've picked to use as the basis for my theory.
    
\section{Does the fact that $H$ depends on $U$ destroy AIT?!}
    
\begin{itemize}
\item
Some people are \textbf{exceedingly unhappy} that the value of $H$ depends
on the choice of $U$.  
But whatever $U$ you pick, the theory goes through.
And Day IV, I actually pick a very simple $U$. 
\item
And since $\Omega$ is an infinite object,
all the additive constants that reflect the choice of $U$ \textbf{wash out in the limit}.
\item
Time complexity depends \textbf{polynomially} on the computational model.
But in the case of program-size complexity the dependence on $U$
is \textbf{only an additive constant}. So the unfashionable field of AIT
depends \textbf{less} on the choice of computer than the fashionable
field of time complexity does!
\textbf{So much for fashion!}
\item
Remember, I do not think that AIT is a practical theory for practical
applications. \textbf{All theories are lies that help us to see the truth.\footnote
{``Art is a lie that helps us to see the truth,'' \emph{Pablo Picasso.}}
All elegant theories are \emph{simplified models} of the horrendous complexity
of the chaotic real world.} Nevertheless, they give us insight.  You have
a choice: elegance or usefulness.  You cannot have both!
\item
If we knew all the laws of physics, and the world turned out to be a giant
computer, then we could \textbf{use that computer} as $U$!  
\end{itemize}

\section{A More Abstract Definition of the Complexity $H$}
    
I prefer a very down-to-earth concrete approach.  
But here is a more abstract ``axiomatic'' approach to defining $H$.
    
Define an abstract complexity measure $H$ via these two properties:
\begin{itemize}
\item
\textbf{$H(x)$ is computable in the limit from above.}
\vspace{5mm}
\\   
I.e., $\{ \langle x,k \rangle : H(x) \leq k \}$ is r.e.
\vspace{5mm}
\\   
I.e., $H(x) = \lim_{t\rightarrow\infty} H^t(x)$.
\item
$\sum_x 2^{-H(x)} \leq 1$.
\end{itemize}

Then pick out an optimal minimal $H$, 
one with the property that for any other 
$H'$ there is a constant $c$ such that
\[  
H(x) \leq H'(x) + c.
\]
    
This \textbf{abstract} approach works!
However I prefer 
the more \textbf{concrete} approach in which 
you 
think of $H(x)$ as the size of the smallest program for $x$.
This gives more \textbf{insight} into what is going on.
Beware of \textbf{premature axiomatization}!
And beware of \textbf{excessive abstraction} concealing \textbf{meaning}!

\section{Back to $U$---What do we do with $U$? Definitions!}
    
\begin{itemize}
\item
Individual complexity of an object. \textbf{Information content} of $x$:
\[   
   H(x) \equiv \min_{U(p)=x} |p|. 
\]   
\item
\textbf{Joint complexity} of a pair of objects $\langle x,y\rangle$: 
\[   
   H(x,y) \equiv \min_{U(p)=\langle x,y\rangle} |p|. 
\]  
Size of the smallest program that computes \textbf{both} objects.
\item
\textbf{Subadditivity}:
\[  
   H(x,y) \leq H(x) + H(y) + c.
\]   
This was the original reason that I made programs self-delimiting, 
so that information content would be additive!\footnote
{\textbf{Then} I discovered $\Omega$, which cannot be defined
unless programs are self-delimiting.}
This means \textbf{we can combine programs} for $x$ and $y$ and 
add $c$ bits and get a program for the pair $\langle x,y\rangle$.
For a LISP program that does this and shows that $c = 432$ is
possible, see Section 4.10.
\item
\textbf{Algorithmic independence} of two objects:
\[  
   H(x,y) \approx H(x) + H(y).
\]  
I.e., the joint complexity is approximately equal to the sum of the
individual complexities.
This means that the two objects have \textbf{nothing in common}.
\end{itemize}

\section{AIT Definitions (Continued)}
    
\begin{itemize}
\item
\textbf{Mutual information} content of two objects:
\[ 
   H(x : y) \equiv H(x) + H(y) - H(x, y).
\]  
This is the extent to which it is better 
to \textbf{compute them together rather than separately}.
\item
An \textbf{elegant program} $y*$ for $y$ is one that is as small as possible. 
I.e., $y*$ has the property that
\[  
   U(y*) = y, \;\; H(y) = |{y*}|, \;\; H(U(y*)) = |{y*}|.
\]  
\item
\textbf{Relative information} content $H(x | y)$ of $x$ 
given an elegant program $y*$ for $y$, \textbf{not $y$ directly}!
\[ 
   H(x | y) \equiv \min_{U(p,y*)=x} |p|. 
\] 
This is the \textbf{second main technical idea of AIT}.  
First is to use self-delimiting programs.
Second is to define relative complexity in this more subtle manner
in which we are \textbf{given $y*$ for free}, not $y$, 
and we have to compute $x$.
\end{itemize}

\section{Important Properties of These Complexity Measures}
    
\begin{itemize}
\item
Here is an \textbf{immediate consequence} of our definitions:  
\[
   H(x, y) \leq H(x) + H(y | x) + c.
\]  
\item
And here is an \textbf{important and subtle theorem} that justifies our definitions!
\[ 
   H(x, y) = H(x) + H(y | x) + O(1).
\] 
This is an identity in Shannon information theory, with no $O(1)$ term.
In these two different versions of information theory,
the formulas look similar but of course the meaning of $H$ is completely different.
\item
Let's apply this subtle theorem to the mutual information:
\[  
   H(x : y) \equiv H(x) + H(y) - H(x, y).
\]  

We get these two corollaries:
\[  
   H(x : y) = H(x) - H(x | y) + O(1),
\]
\[
   H(x : y) = H(y) - H(y | x) + O(1).
\]
    
So the \textbf{mutual information is also \emph{the extent to which knowing
one of the two objects helps you to know the other}}.  
It was not at all
obvious that this would be symmetric!
    
Again, these two corollaries are identities 
in Shannon information
theory.
They justify our whole approach.
When I got these results \textbf{I knew 
that the definitions I had picked
for AIT were finally correct}!
\end{itemize}

\section{Complexity of Bit Strings. Examples.}

\begin{itemize}
\item
\textbf{High complexity:}
\begin{eqnarray*}
   H(\mbox{\textbf{the most complex $N$-bit strings}}) 
     & = & N + H(N) + O(1) 
\\
     & \approx & N + \log_2 N.
\end{eqnarray*}
In other words, you have to know how many bits there are
plus what each bit is.
And most $N$-bit strings have $H$ close
to the maximum possible.
These are the \textbf{algorithmically random} $N$-bit strings.

This notion has many \textbf{philosophical resonances}.  First, a random
string is one for which there is no theory that obeys Leibniz's dictum that
a theory must be smaller than what it explains.  
\textbf{Leibniz clearly anticipated this definition of randomness.} 
   
Second,
a random string is an unexplainable, irrational string, one that cannot be comprehended,
except as ``a thing in itself'' (\textbf{Ding an sich}), to use Kantian
terminology.
\item
\textbf{Low complexity:}
\\
Consider the $N$-bit strings consisting of $N$ 0's or $N$ 1's.
\[  
   H(0^N) \;\; = \;\; H(N) + O(1) \;\; \approx \;\; \log_2 N,
\]  
\[
   H(1^N) \;\; = \;\; H(N) + O(1) \;\; \approx \;\; \log_2 N.
\]  
For such bit strings, 
you only need to know how many bits there are,
not what each bit is.
These bit strings are \textbf{not at all} random.
\item
\textbf{Intermediate complexity:}
\\
Consider an elegant program $p$.  
I.e., no smaller program makes $U$ produce the same output.
\[
   H(p) = |p| + O(1) \;\;\; \mbox{\textbf{if $p$ is an elegant program}.}
\]
These bit strings are \textbf{borderline} random.
Randomness is a matter of degree,
and this is a good place to put the cut-off 
between random
and non-random strings, 
if you have to pick a cut-off.
\end{itemize}

\textbf{How about \emph{infinite} bit strings?}

\section{Maximum Complexity Infinite Binary Sequences}
    
\begin{itemize}
\item
\textbf{Maximum Complexity:}
\\  
An infinite binary sequence $x$ 
is defined to be \emph{algorithmically
random} iff
\[  
   \exists c \; \forall N [ \; H(x_N) > N - c \; ].
\] 
Here $x_N$ is the first $N$ bits of $x$.
   
If $x$ is generated by independent tosses of a fair coin, 
then with probability one, $x$ is algorithmically random.
But how about specific examples?
\item 
\textbf{Two Examples of Maximum Complexity:}
\[  
\Omega \equiv \sum_{\mbox{$U(p)$ halts}} 2^{-|p|} 
\]
\[
\Omega' \equiv \sum_N 2^{-H(N)} 
\]
These are algorithmically random real numbers.
The initial segments of their binary expansion
always have \textbf{high} complexity.
    
From this it follows that $\Omega$ and $\Omega'$ aren't
contained in any constructively definable set of measure zero.
In other words, they are statistically (Martin-L\"of) random.
Therefore
$\Omega$ and $\Omega'$ 
\textbf{\emph{necessarily have}} 
any property held by infinite binary
sequences $x$ with probability one.
E.g., $\Omega$ and $\Omega'$ are both \textbf{Borel normal},
\textbf{\emph{for sure}},
since in general this is the case with probability one.
Borel normal means that in any base, in the limit all blocks of digits of the
same size provably have equal relative frequency.

That $\Omega$ and $\Omega'$ are both Borel normal
real numbers
can also be proved directly using a simple program-size argument.
This can be done
because, as Shannon information theory already shows, 
reals that
are not Borel normal are highly redundant and compressible.
And the main difference between Shannon information theory and AIT
is that Shannon's theory is concerned with \textbf{average} compressions
while I am concerned with compressing \textbf{individual} sequences.
\item
\textbf{Why can't we demand the complexity of initial segments
be even higher?}
\\  
Because with probability one, 
infinitely often
there have to be \textbf{long runs} 
of consecutive zeroes or ones,
and this makes the complexity dip far below 
the maximum
possible [$N + H(N)$] infinitely often! 
\end{itemize}

\section{The Algorithmic Probability of $x$, $P(x)$}
    
Straight-forward theorems in AIT use program-size arguments.
    
But the more subtle proofs have to descend one level and use probabilistic arguments.
AIT is an iceberg!  Above the water we see program size.
But the bulk of the AIT iceberg is submerged and is the algorithmic probability $P(x)$!
   
$P_C(x)$ is defined to be 
the probability that a program 
generated by coin-tossing
makes $C$ produce $x$:
\[  
   P_C(x) \equiv \sum_{C(p)=x} 2^{-|p|},
\]
\[
   P(x) \equiv P_U(x) \equiv \sum_{U(p)=x} 2^{-|p|}.
\]  
This takes into account \textbf{all} the programs that calculate $x$,
not just the smallest one.

\section{What is the Relationship Between $H$ and $P?$}
    
First of all, it's obvious that
\[ 
   P(x) \geq 2^{-H(x)}
\] 
because one way to compute $x$ is using an elegant program for $x$.
In fact, much, much more is true:
\[ 
   H(x) = -\log_2 P(x) + O(1).
\]
This crucial result shows that AIT, 
at least in so far as the appearance of its \textbf{formulas} is concerned,
is sort of just a version of probability theory
in which we take logarithms and convert probabilities into complexities.
In particular, the important and subtle theorem that
\[ 
   H(x, y) = H(x) + H(y | x) + O(1)
\]
is seen from this perspective as merely 
an alternative version of the \textbf{definition} of relative probability:
\[
   \mathrm{Pr}(x, y) \equiv \mathrm{Pr}(x) \times \mathrm{Pr}(y | x).
\]  
How do we establish this crucial relationship between $H$ and $P$?
By using an extended version of something called the Kraft inequality,
which AIT inherited from Shannon information theory!
[See Section 3.13.]

\section{Occam's Razor!  There are few elegant programs!}
    
But first I want to discuss an important corollary of the crucial
theorem 
\[ 
   H = -\log_2 P
\] 
connecting $H(x)$ and $P(x)$.
    
This shows that there cannot be too many
small programs for $x$. 
    
In particular, it follows that
the number of elegant programs for $x$ is bounded.
    
Also,
the number of programs for $x$ 
whose size is within $c$ bits of $H(x)$ 
is bounded by a function that depends only on $c$ and not on $x$.
   
This function of $c$ is approximately $2^c$.
    
In fact,
any program for $x$ whose size is close to $H(x)$ 
can
easily be obtained from any other such program.
    
So, what I call Occam's razor, elegant or nearly elegant programs
for $x$ \textbf{are essentially unique}!
    
And this justifies our definition of relative complexity $H(x | y)$
and shows that it does not depend too much on the choice
of the elegant program for $y$.
There are at most $O(1)$ bits difference in
$H(x | y)$
depending on the choice of the $y*$ that we are given for free.

\section{The Extended Kraft Inequality Condition for Constructing $C$}
    
The crucial tool used to show that
\[ 
   H(x) = -\log_2 P(x) + O(1)
   \;\;\;\; \mbox{and} \;\;\;\;
   H(x, y) = H(x) + H(y | x) + O(1)
\]  
is an extended version of the Kraft inequality condition for the existence
of a prefix-free set of words.   
    
In AIT, the Kraft inequality gives us a necessary and sufficient condition not for the
existence of a prefix-free set, 
as it did in Shannon information theory,
but for the existence of a self-delimiting computer $C$.
Once we have constructed a special-purpose computer $C$ using the Kraft
inequality, we can then make statements about $U$ by using the fact that
\[ 
   H_U \leq H_C + c.
\] 

Here's how we construct $C$!
    
Imagine that we have an algorithm for generating 
an infinite list of
requirements:
\[ 
   \langle\mbox{size of program for $C$}, \mbox{output that we want from $C$}\rangle
\] 
As we generate each requirement $\langle s, o\rangle$,
we pick \textbf{the first available} $s$-bit program $p$ 
and we assign the output $o$ to $C(p)$. 
Available means not an extension or prefix of 
any previously assigned program for $C$.
This process will work and will produce a self-delimiting 
computer $C$ satisfying all the 
requirements iff 
\[ 
   \sum_{\mbox{over all requirements $\langle s, o\rangle$}} 2^{-s} \leq 1.
\]  
Note that if there are duplicate requirements in the list, 
then several $p$ will yield the same
output $o$. 
I.e., we will have several $p$ with the same value for $C(p)$.
    
This way of constructing $C$ may be thought of as a first-fit storage 
allocation
algorithm for one infinitely-divisible unit of storage.

\section{Three-Level Proofs: Program Size, Probabilities, Geometry}
   
The proof of this crucial version of the Kraft inequality 
involves simultaneously thinking of $C$ as a computer, as an assignment
of probabilities to outputs, and as an assignment of outputs to
segments of the unit interval, segments which are halves, or halves of halves,
or halves of halves of halves\ldots
    
It is important to be able to simultaneously keep each of these three images in ones mind,
and to translate from one image to the next depending on which gives the most
insight at any given moment.
    
In other words, $C$ may be thought of as a kind of constructive probability
distribution, and probability one corresponds to the entire unit interval
of programs, with longer and longer programs corresponding to smaller
and smaller subintervals of the unit interval.
    
Then the crucial fact that no extension of a valid program is a valid program
simply says that the intervals corresponding to valid programs
\textbf{do not overlap}.

In the last lecture, Day IV, I'll drop this abstract viewpoint and make everything
very, very concrete by indicating how to program $U$ in LISP and actually run
programs on $U$\ldots

\chapter{Day IV---LISP Implementation of AIT}
\markboth{CHAITIN, FROM PHILOSOPHY TO PROGRAM SIZE}{DAY IV---LISP IMPLEMENTATION}

The LISP interpreter for Day IV, 
which is a Java applet that will run in your web browser,
can be found at this URL:
\vspace{4mm}
\\
\texttt{http://www.cs.auckland.ac.nz/CDMTCS/chaitin/unknowable/lisp.html}
\vspace{0mm}
\\
Note that there is a concise LISP ``reference manual'' on these web pages,
just below the windows for input to and output from the interpreter.
     
\section{Formal Definition of the Complexity of a Formal Theory}

Before starting with LISP, let me finish some things
I started in the last lecture, and also answer some
questions that were raised about that lecture.
    
First of all, I never stated the formal version of
the incompleteness results about $\Omega$ that
I explained informally in Day II.
    
\textbf{How can we measure the complexity of a formal
mathematical theory in bits of information?}
    
$H(\mbox{formal math theory}) \equiv$ 
\[
   \mbox{the size in bits $|p|$ of the smallest (self-delimiting!)\@ program $p$} 
\]
such that
\[ 
   U(p) = \{\mbox{the infinite set of all the theorems in the theory}\}.
\]
    
Note that this is something new: $U$ is now performing
an \textbf{unending} computation and produces \textbf{an infinite 
amount} of output!
\[               
   \mbox{smallest $p$ $\longrightarrow$ $U$ $\longrightarrow$  \{theorems\}}
\]

The rather abstract point of view taken here is in
the spirit of Emil Post's 1944 paper \emph{R.e.\ sets of
positive integers and their decision problems.}
For our purposes the internal details of the formal 
theory are irrelevant.  We are flying high in the sky
looking down at the formal theory. We're 
so high up that we don't care about the axioms
and rules of inference used in the theory. 
We can't see all of that detail. Nevertheless,
using these very general methods we can make some
rather strong statements about the power of the formal
theory in question.

\section{Formal Statement of the Incompleteness Theorem for $\Omega$}

\textbf{A formal math theory $T$ with complexity $N$
cannot enable you to determine more
than $N + c$ bits of the numerical value of
$\Omega$ (base-two).}
    
Here are some hints about the proof.
We make the tacit hypothesis, of course,
that 
all the theorems proved in theory $T$ are true.
Otherwise $T$ is of no interest.
    
The proof is in two parts.
\begin{itemize}
\item
First of all you show that
\[  
H(\Omega_K) > K - \mbox{constant}.
\]  
This follows from the fact that
knowing $\Omega_K$,
the first $K$ bits of $\Omega$,
enables you to answer the halting
problem for all programs for $U$
up to $K$ bits in size.
\item
Now that you know that
$\Omega_K$
has high complexity, in
the second half of the proof you
use this fact to derive 
a contradiction
from the assumption that the formal
theory $T$ enables you to determine
much more than 
$H(T)$ bits of  
$\Omega$.
Actually, it makes no difference
whether these are consecutive bits
of $\Omega$ or are scattered 
about
in the base-two expansion of $\Omega$.
\end{itemize}
    
This, our major incompleteness
result about $\Omega$, 
actually requires
very little of the machinery presented
in Day III.  
In fact, the proof is a
straight-forward program-size argument
that makes no use of $P$ nor of the
the Kraft inequality.  
This is in fact the self-contained
elementary
proof that I give 
in my 1998 Springer volume
\emph{The Limits of Mathematics},
which is the reference volume for today's lecture.

\section{Mutual Information is Program Size?}
    
I was asked a question about mutual information.
    
Day III we defined the following complexity measures:
\begin{itemize}
\item
$H(x)$ (\emph{absolute, individual}) information,
\item
$H(x, y)$ \textbf{joint} information,
\item
$H(x | y)$ \textbf{relative} information,
\item
$H(x : y)$ \textbf{mutual} information.
\end{itemize}
The fourth complexity measure in this list, mutual information,
is actually rather different from the other three.  
\textbf{It isn't the size of a computer program}!
Instead it's defined to be
\[ 
   H(x) + H(y) - H(x, y).
\]   

Well, you \emph{could} try to make mutual information into the size of
a program.
E.g., mutual information could be defined to be the size of
the \textbf{largest} (not the smallest!)\@ possible common subroutine
shared by elegant programs for $x$ and $y$.\footnote
{For other possibilites, see the section on \textbf{common information}
in my 1979 paper \emph{Toward a mathematical definition of ``life.''}}
    
Unfortunately this alternative definition, which is much more intuitive
than the one I actually use, doesn't seem to fit into AIT too well.
Perhaps someone can find a way to make it work!  Can you?
    
This is an interesting question to answer because of an important
possible application: In Section 1.14, Day I, it was suggested that
mutual information be used to try to partition the world into 
separate living organisms.

\section{How do you actually make things self-delimiting?!}

Is 
\[
   H(\mbox{$N$-bit string}) \leq N + \mbox{constant}? 
\]
No, it can't be, because programs for $U$ have to be self-delimiting.
\textbf{Programs have to indicate their size as well as their content.}
     
Let's look at some examples of how this can be done.
How can we make an $N$-bit string self-delimiting?
Well, $2N + 2$ bits will certainly do.
There's a special-purpose 
self-delimiting computer $C$ that accomplishes this:
\[ 
   H_C(\mbox{$N$-bit string}) \leq 2N + 2.
\]  
$C$ reads two bits of its program 
at a time as long as they are the same,
then stops reading its program 
and outputs one bit from 
each pair of identical twins that it read.
So
\[ 
 C(00\; 00\; 11\; 01) \longrightarrow 001.
\] 
Therefore
\[ 
   H(\mbox{$N$-bit string}) \equiv
   H_U(\mbox{$N$-bit string}) \leq 
   H_C(\mbox{$N$-bit string}) + c \leq 
   2N + c'.
\] 
A more complicated $C$ reads a program that starts with a doubled
and therefore self-delimiting version of the base-two numeral
for $N$ followed by the $N$ bits of the bit string. This shows that
\[ 
   H(\mbox{$N$-bit string}) \leq 2\log_2 N + N + c.
\]  
A more elaborate procedure is to have \textbf{two} ``headers'' at the beginning of the
program.  The first one uses the bit-doubling trick to indicate the size of the
second header, and the second header gives $N$ in binary directly, with
no bits doubled.  That gives an upper bound of 
\[  
   H(\mbox{$N$-bit string}) \leq O(\log\!\log N) + \log_2 N + N.
\] 
You can go on like this using more and more headers and getting more and
more complicated upper bounds.  To short-circuit this process, just go back
to using a \textbf{single} header, and make it an elegant program for $N$. 
That gives
\[  
   H(\mbox{$N$-bit string}) \leq H(N) + N + c,
\] 
and this upper bound is, in general, best possible.
    
So we've just seen four different examples of how a self-delimiting 
program can tell us its
size as well as its contents.
    
Another version of the principle that a program tells
us its size as well as its contents is the fact that 
an elegant program tells us its size as well as its output, so
\[ 
   H(x) = H(x, H(x)) + O(1).
\]  
\emph{This is one formula in AIT that is totally unlike anything 
in Shannon information theory.}
   
\textbf{That's all the unfinished business from Day III!
Now let's start Day IV proper by explaining LISP!}

\section{What is LISP?}
     
LISP = \emph{LISt Processing.} Lists are arbitrarily long nested \textbf{tuples},
and may have repeated elements.
    
LISP is 
like a computerized version of set theory!
     
The LISP S-expression 
\begin{center}
   \texttt{((A bc) 39 x-y-z)} 
\end{center}
denotes the following nested tuples:
\begin{center}
   $\langle\langle \mathtt{A},\;\mathtt{bc}\rangle,\;\mathtt{39},\;\mbox{\texttt{x-y-z}}\rangle$
\end{center}
Programs and data in LISP are always S-expressions.
Everything is an S-expression!
The S in S-expression stands for \emph{Symbolic.}
    
S-expressions that are not \textbf{lists} are called \textbf{atoms}.
\texttt{A}, \texttt{bc}, \texttt{39} and \texttt{x-y-z} 
are the \emph{atoms} in the above S-expression.
The empty list \texttt{()} = $\langle\rangle$ is also an atom, usually
denoted by \texttt{nil}.

LISP is a functional programming language,
not an imperative programming language.
    
LISP is an expression language. There is no notion of time,
there are no \textbf{GOTO}'s, and there are no assignment statements!\footnote
{Well, they're actually present, but only indirectly, only subliminally!}
    
LISP programs are expressions that you evaluate,
not run, and they yield a value, with \textbf{NO} side-effect!
(This is called ``Pure'' LISP!)
    
\newpage
\textbf{Functional Notation in LISP:}
\begin{center}
$f(x, y)$ is written \texttt{(f x y)}
\end{center}
So 
\begin{center}
   $(1 + 2) \times 3$ 
\end{center}
becomes
\begin{center}
   \texttt{(* (+ 1 2) 3)} 
\end{center}
in LISP.
Prefix notation, not infix notation! Full parenthesization!
Fixed number of arguments for each primitive function!
    
In my LISP most parentheses are understood, are implicit,
and are then supplied by the LISP interpreter.
     
\textbf{Complete LISP Example.} 
Factorial programmed in my LISP:
{\large
\begin{verbatim}
   define (f n)
   if = n 0  1
      * n (f - n 1)
\end{verbatim}
}
     
With all the parentheses supplied, factorial becomes:
{\large
\begin{verbatim}
   (define (f n)
   (if (= n 0)  1
       (* n (f (- n 1)))
   ))
\end{verbatim}
}

In words, $f$ is defined to be a function of one argument, $n$.
And if $n$ = 0, then $f(n)$ is 1. 
If not, $f(n)$ is 
$n\times f(n-1)$.
     
Then \texttt{(f 4)} yields value 24, as expected, since  
$4\times 3\times 2\times 1 = 24$.

\textbf{Dissecting and Reassembling Lists:}
\\
Finally I should mention that \texttt{car} gives you the first element of a list,
\texttt{cdr} gives you the rest of the list, and \texttt{cons} inverts \texttt{car} and \texttt{cdr}.

\section{LISP Program-Size Complexity}
    
You can get a toy version of AIT by using LISP program size!
    
Define an elegant LISP expression to be one with
the property that no smaller LISP expression
yields the same value.
    
In order to measure size properly,
LISP expressions are written in a canonical standard
notation with no embedded comments and with exactly one blank
separating successive elements of a list. 
    
Then the \textbf{LISP complexity} of an S-expression is defined
to be the size in characters of an elegant expression
with that particular value.
    
Next represent a formal mathematical theory in LISP
as a LISP S-expression whose evaluation never terminates
and which uses the primitive function \textbf{display} to 
output each theorem.
    
\textbf{display} is an identity function with
the side-effect of displaying its operand and is normally 
used for debugging large LISP S-expressions
that give the wrong value.
    
\textbf{\emph{Theorem---}
\\
A formal mathematical theory with LISP complexity $N$ cannot
enable you to prove that a LISP expression is elegant
if the elegant expression's size is greater than $N + 410$ characters!}
    
This is a very concrete, down-to-earth incompleteness theorem!
In that direction, it's the best I can do!
    
The \emph{reductio ad absurdum}
proof consists of $N + 410$ characters of LISP that define a function,
and that apply this function to a quoted form of the $N$-character LISP expression
that generates all the theorems of the formal theory.  The first
step is to measure the size $N$ of the LISP expression that generates
all the theorems.  Then use TRY (Section 4.7) to run the theory for longer
and longer, until a provably elegant expression is found whose size
is \textbf{greater} than $N + 410$. When this happens, evaluate the provably
elegant expression and return its value as our value. 
    
So we've described an $N + 410$ character LISP expression whose value is the
same as the value of an elegant expression that is larger than it is!
Contradiction!

\section{How Does TRY work?}
    
TRY plays the fundamental role in my LISP that 
\texttt{eval}
plays in traditional, normal LISP. It enables you
to try to evaluate a given expression for a limited
or unlimited amount of time, and giving it raw binary
data on the side, which it must access in a self-delimiting manner.
Also, 
\texttt{display}'s 
from within the expression are captured,
thus giving a mechanism for an S-expression to produce
an infinite amount of output, not just a final value.
     
\begin{verbatim}
   try time-limit/no-time-limit expression binary-data
\end{verbatim}
     
\textbf{TRY has the three arguments shown above.}
     
The first argument is either an integer or 
\texttt{no-time-limit}.
The second argument is an arbitrary LISP S-expression.
The third argument is a bit string represented in LISP
as a list \texttt{(...)} of 0's and 1's separated by blanks.

\textbf{TRY then yields this triple:}
\begin{verbatim}
   (success/failure
    value-of-expression/out-of-time/out-of-data
    (list of captured displays from within expression...))
\end{verbatim}

Also, within the expression that is being tried,
you can use 
\texttt{read-bit}
or 
\texttt{read-exp}
to get access
to the binary data in a self-delimiting manner,
either one bit or one LISP S-expression at a time.
    
TRY is like a large suitcase: I am using it to do
several different things at the same time. In fact,
TRY is all I really need to be able to program AIT in
LISP.  All of the other changes in my version of LISP
were made just for the fun of it! (Actually, to
\textbf{simplify} LISP as much as possible without
ruining its power, so that I could prove theorems
about it and at the same time enjoy programming in it!)
     
\textbf{Example:} If you're trying to run an unending
computation that is a formal theory written in LISP
as described in Section 4.6, then TRY always returns
\begin{verbatim}
   (failure out-of-time (list of theorems...))
\end{verbatim}
The size of the list of theorems depends on how much
time the TRY was given to run the formal theory.

\section{LISP Program for Our Standard Universal Computer $U$}

{\large
\begin{verbatim}
   define (U p)
   cadr try no-time-limit
            'eval read-exp
            p
\end{verbatim}
}
     
This function $U$ of $p$ returns the second element of the list returned by TRY,
a TRY with no time limit that tries to evaluate a LISP expression at the
beginning of $p$, while making the rest of $p$ available as raw binary data.
    
Here the program $p$ for $U$ is a bit string 
that is
represented 
in LISP as a list of 0's and 1's separated by blanks.
E.g., $p$ = \texttt{(0 1 1 0 1...)}
    
Let me explain this!
    
The binary program $p$ consists of a LISP S-expression prefix (converted into 
a bit string), that is followed by a special delimiter character, the ASCII newline
control character (also converted to bits), and that is then followed by raw binary data.
    
The newline character adds 8 bits to the program $p$, but guarantees
that the prefix be self-delimiting.  Each character in the LISP prefix is
given in $p$ as the corresponding 8 ASCII bits.
    
Within the prefix, \texttt{read-bit} and 
\texttt{read-exp} enable you to get access to
the raw binary data, either one bit or one LISP S-expression at a time.  
But you must not run off the end of the binary data!
That aborts everything! If so, $U$ of $p$ returns \texttt{out-of-data}.
    
So the program for the universal $U$ consists of a high-level language algorithmic
part followed by data, raw binary data.  The algorithmic part indicates
which special-purpose self-delimiting computer $C$ should be simulated,
and the raw binary data gives the binary program for $C$! 

\textbf{I believe that I've picked a \emph{very natural} $U$.}

\section{Running $U$ on Simple Examples}
    
Here's our first example of a program for $U$.
This program is all LISP prefix with no binary data.
The prefix is \texttt{'(a b c)}.
\begin{verbatim}
   run-utm-on
   bits ' '(a b c)
\end{verbatim}
This yields 
\texttt{(a b c)}.
     
And here's a program with one bit of binary data.
The prefix is \texttt{read-bit} and the data is \texttt{0}.
\begin{verbatim}
   run-utm-on
   append bits 'read-bit 
          '(0)
\end{verbatim}
This yields \texttt{0}.
   
Let's change the data bit.
Now the prefix is \texttt{read-bit} and the data is \texttt{1}.
\begin{verbatim}
   run-utm-on
   append bits 'read-bit 
          '(1)
\end{verbatim}
This yields \texttt{1}.
    
\texttt{run-utm-on} is a macro that expands to the definition of $U$,
which is \texttt{cadr try no-time-limit 'eval read-exp}.
    
Also note that when \texttt{bits} converts a LISP S-expression
into a bit string (a list of 0's and 1's) it automatically
adds the newline character required by \texttt{read-exp}.

\section{Subadditivity Revisited}
     
{\large
\begin{verbatim}
   cons eval read-exp
   cons eval read-exp
        nil
\end{verbatim}
}
This is a 432-bit prefix $\pi$ for $U$
with the property that 
\[  
   U(\pi \, {x*} \, {y*}) = (x \; y) = \langle x, y\rangle.
\]  
Here $x*$ is an elegant program for $x$ and
$y*$ is an elegant program for $y$.
    
\textbf{Therefore}
\[   
   H(x, y) \leq H(x) + H(y) + 432 !
\]  

For years and years I used this inequality
without having the faintest idea what the value of
the constant might be!  It's nice to know that
there is a natural choice for $U$ that can be easily
programmed in LISP and for which $c$ = 432!  That's
a little bit like dreaming of a woman for years
and then actually meeting her! That's also how I feel
about now having a version of $U$ that I can actually
run code on!
    
My book \emph{The Limits of Mathematics} programs in LISP
all the proofs of my incompleteness results, in particular,
all the key results about $\Omega$.  And my book
\emph{Exploring Randomness} programs in LISP all the proofs
of the results presented Day III.  That's a lot of 
LISP programming, but I enjoyed it a lot!  

\textbf{Some} of
these programs are fun to read, and it's educational to do so,
but I probably went overboard!  In general, I would say that
you are better off programming something yourself, rather
than reading someone else's code, if you \textbf{really} want
to understand what's going on!  There is no substitute for
doing it yourself, your own way, if you really want to 
understand something.  After all, that's how I came up
with AIT and my alternative approach to G\"odel incompleteness
in the first place, because I was dissatisfied with the usual
approach!

\section{Some Open Problems for Future Research}

\begin{itemize}
\item
There is a cottage industry of highly-technical work
on subtle properties of $\Omega$ and to what 
extent other numbers share these properties.
[See 
Delahaye, ``Les nombres om\'ega,'' \emph{Pour la Science}, May 2002, pp.\ 98--103 
and
Calude, \emph{Information and Randomness}, Springer, 2002.]
\item
Much work remains to be done on the complexity $H(X)$ for
\textbf{infinite} sets $X$. 
This is important because of the case when
\[
   X = \{\mbox{the set of theorems of a formal theory $T$}\}.
\]
For example, the relationship 
between the algorithmic complexity $H(X)$ and 
the algorithmic probability
$P(X)$ is not completely known.
When dealing with infinite computations in AIT,
in some cases it is not even clear whether the \textbf{definitions} 
that were chosen are the right ones!
[For more on this, see the last chapter of my \emph{Exploring Randomness.}]
\item
As I stated above, meta-mathematics and the incompleteness results
I've presented in this course
are a \textbf{static framework} in which one studies the power of a \textbf{fixed} set
of axioms and rules of inference. How about a \textbf{dynamic} theory that has
something to say about the \textbf{emergence of new mathematical concepts}?
Is such a meta-theory possible?\footnote
{In this connection, see footnote 6 on page 24.}
\item
I would like to see an abstract mathematical theory defining ``life''
and predicting its emergence and \textbf{evolution} under very general circumstances.
\item
Is the physical universe \textbf{discrete} or \textbf{continuous}? Starting from
black hole thermodynamics, a ``holographic'' principle has emerged which
states that the information content of any physical system is only a finite
number of bits 
and grows with the \textbf{surface area} of that system, not with its volume.
[See Smolin, \emph{Three Roads to Quantum Gravity.}]
It will be interesting to see if and in what form these ideas survive.
\item
I think that AIT is largely sketched out and only highly technical questions
remain.  
I personally am more interested in attempting \textbf{new theories} of life or intelligence.
This however will not be easy to do.
\end{itemize}
\textbf{I've greatly enjoyed this Winter School and my visit to Estonia!
Thank you very much for inviting me!}

\chapter*{Additional Reading}
\addcontentsline{toc}{chapter}{Additional Reading}
\markboth{CHAITIN, FROM PHILOSOPHY TO PROGRAM SIZE}{ADDITIONAL READING}

The six papers listed here were the handouts that accompanied my lectures at EWSCS '03.
The two books listed below provide the reference material for the course,
including complete proofs for all the theorems that I stated here,
and a detailed explanation of my version of LISP.

\begin{itemize}
\item
\textbf{Day I---Philosophical Necessity of AIT}
\begin{itemize}
\item
G. Chaitin, ``Two philosophical applications of algorithmic information theory,''
\emph{Proc.\ of DMTCS '03}, Springer-Verlag, to appear.
\\
\texttt{http://arxiv.org/abs/math.HO/0302333}
\item
G. Chaitin, ``On the intelligibility of the universe and the notions of simplicity,
complexity and irreducibility.''
\\
Presented at the XIX.\ Deutscher Kongre{\ss} f\"ur
Philosophie, \emph{Grenzen und Grenz\"uberschreitungen}
[German Philosophy Congress on Limits and Transcending Limits], 23.--27.\ September 2002 in Bonn.
\\
\texttt{http://arxiv.org/abs/math.HO/0210035}
\end{itemize}
\item
\textbf{Day II---Main Application of AIT: Incompleteness}
\begin{itemize}
\item
G. Chaitin, ``Paradoxes of randomness,'' \emph{Complexity}, vol.\ 7, no.\ 5, pp.\ 14--21, 2002.
\end{itemize}
\item
\textbf{Day III---Technical Survey of AIT: Definitions \& Theorems}
\begin{itemize}
\item
G. Chaitin, ``Meta-mathematics and the foundations of mathematics,'' \emph{Bulletin EATCS},
vol.\ 77, pp.\ 167--179, 2002.
\item
G. Chaitin, \emph{Exploring Randomness,} Springer-Verlag, London, 2001.
\end{itemize}
\newpage
\item
\textbf{Day IV---LISP Implementation of AIT}
\begin{itemize}
\item
G. Chaitin, ``Elegant LISP programs,'' in C. S. Calude, ed., \emph{People and Ideas in
Theoretical Computer Science}, pp.\ 32--52. Springer-Verlag, Singapore, 1999.
\item
G. Chaitin, ``An invitation to algorithmic information theory,'' in D. S. Bridges, C. S.
Calude, J. Gibbons, S. Reeves, I. H. Witten, eds.,  \emph{Combinatorics, Complexity
\& Logic: Proc.\ of DMTCS '96}, pp.\ 1--23. Springer-Verlag, Singapore, 1997.
\item{}
[The above two papers are both chapters in]
\\
G. Chaitin, \emph{The Limits of Mathematics,} Springer-Verlag, Singapore, 1998.
\end{itemize}
\end{itemize}
The papers listed here are also available at the author's website:
\vspace{2mm}
\\
\texttt{http://www.cs.auckland.ac.nz/CDMTCS/chaitin}
\vspace{2mm}
\\
LISP code for my two books may be found at
\vspace{2mm}
\\
\texttt{http://www.cs.auckland.ac.nz/CDMTCS/chaitin/ait}
\\
\texttt{http://www.cs.auckland.ac.nz/CDMTCS/chaitin/lm.html}
\vspace{2mm}
\\
The LISP interpreter is at
\vspace{2mm}
\\
\texttt{http://www.cs.auckland.ac.nz/CDMTCS/chaitin/unknowable/lisp.html}

\end{document}